\documentclass[12pt]{amsart}
\usepackage{enumerate}
\usepackage{amsmath, amsfonts, amssymb, amsthm}
 \makeatletter
\@namedef{subjclassname@2010}{%
  \textup{2010} Mathematics Subject Classification}
\makeatother
\newtheorem{thm}{Theorem}[section]

\newtheorem{lem}[thm]{Lemma}

\newtheorem{conj}[thm]{Conjecture}

\theoremstyle{definition}

\newtheorem*{xrem}{Remark}

\numberwithin{equation}{section}

\def\Z{\Bbb Z}
\def\N{\Bbb N}

\def\l{\left}
\def\r{\right}
\def\bg{\bigg}
\def\({\bg(}
\def\){\bg)}
\def\]{]\!]}
\def\[{[\![}
\def\t{\text}
\def\f{\frac}
\def\mo{{\rm{mod}\ }}

\def\ls{\leqslant}
\def\gs{\geqslant}
\def\se {\subseteq}

\def\eq{\equiv}

\begin{document}

\baselineskip=17pt

\hbox{Colloq. Math. 145(2016), no.\,1, 149-155.}
\medskip

\title
[On some universal sums of generalized polygonal numbers] {On some universal sums of generalized polygonal numbers}

\author
[F. Ge and Z.-W. Sun] {Fan Ge and Zhi-Wei Sun}

\address{(Fan Ge) Department of Mathematics\\ University of Rochester
\\Rochester, NY 14627, U.S.A.} \email{fange.math@gmail.com}
\address {(Zhi-Wei Sun, corresponding author) Department of Mathematics\\ Nanjing
University\\ Nanjing 210093\\ People's Republic of China}
\email{zwsun@nju.edu.cn}

\date{}

\begin{abstract}
 For $m=3,4,\ldots$ those $p_m(x)=(m-2)x(x-1)/2+x$ with
$x\in\Z$ are called generalized $m$-gonal numbers. Sun \cite{S15}
studied for what values of positive integers $a,b,c$
the sum $ap_5+bp_5+cp_5$ is universal over $\Z$ (i.e., any
$n\in\N=\{0,1,2,\ldots\}$ has the form $ap_5(x)+bp_5(y)+cp_5(z)$
with $x,y,z\in\Z$). We prove that
$p_5+bp_5+3p_5\,(b=1,2,3,4,9)$ and $p_5+2p_5+6p_5$ are universal
over $\Z$, as conjectured by Sun. Sun also conjectured that any $n\in\N$ can be
written as $p_3(x)+p_5(y)+p_{11}(z)$ and $3p_3(x)+p_5(y)+p_7(z)$
with $x,y,z\in\N$; in contrast, we show that $p_3+p_5+p_{11}$ and
$3p_3+p_5+p_7$ are universal over $\Z$. Our proofs are essentially elementary and hence suitable for general readers.
\end{abstract}



\maketitle

\section{Introduction}\setcounter{equation}{0}

 For $m=3,4,\ldots$ we set
\begin{equation}\label{1.1}p_m(x)=(m-2)\f{x(x-1)}2+x.\end{equation}
Those $p_m(n)$ with $n\in\N=\{0,1,2,\ldots\}$ are the well-known
{\it $m$-gonal numbers} (or polygonal numbers of order $m$). We call those
$p_m(x)$ with $x\in\Z$ {\it generalized $m$-gonal numbers}. Note that
(generalized) 3-gonal numbers are triangular numbers and
(generalized) 4-gonal numbers are squares of integers.

In 1638, Fermat asserted that each $n\in\N$  can be written as the sum
of $m$ polygonal numbers of order $m$. This was proved by Lagrange,
Gauss and Cauchy in the cases $m=4$, $m=3$ and $m\gs5$ respectively
(see Moreno and Wagstaff \cite[pp.\,54-57]{MW} or Nathanson \cite[Chapter 1, pp.\,3-34]{N96}).
The generalized pentagonal numbers play a crucial role in Euler's famous recurrence for the partition
function.

For $a,b,c\in\Z^+=\{1,2,3,\ldots\}$ and $i,j,k\in\{3,4,\ldots\}$,
Sun \cite{S15} called the sum $ap_i+bp_j+cp_k$ universal over $\N$
(resp., over $\Z$) if for any $n\in\N$ the equation
$n=ap_i(x)+bp_j(y)+cp_k(z)$ has solutions over $\N$ (resp., over
$\Z$). In 1862 Liouville (cf. \cite[p.\,23]{D99}) determined all those
universal $ap_3+bp_3+cp_3$. The second author \cite{S07} initiated the determination
of those universal sums $ap_i+bp_j+cp_k$ with $\{i,j,k\}=\{3,4\}$, and this project was
completed via \cite{S07,GPS,OS}. For almost universal sums $ap_i+bp_j+cp_k$ with $\{i,j,k\}\se\{3,4\}$, see \cite{KS,CH,CO}.

It is known that generalized hexagonal
numbers are identical with triangular numbers (cf. \cite{G94} or \cite[(1.3)]{S15}).

The second author recently established the following result.

\begin{thm}\label{Th1.1} {\rm (Sun \cite[Theorem 1.1]{S15})} Suppose that
$ap_k+bp_k+cp_k$ is universal over $\Z$, where
$k\in\{4,5,7,8,9,\ldots\}$, $a,b,c\in\Z^+$ and $a\ls b\ls c$. Then
$k=5$, $a=1$ and $(b,c)$ is among the following $20$ ordered pairs:
\begin{align*}&(1,c)\ (c\in\{1,2,3,4,5,6,8,9,10\}),
\\&(2,2),\ (2,3),\ (2,4),\ (2,6),\ (2,8),
\\&(3,3),\ (3,4),\ (3,6),\ (3,7),\ (3,8),\ (3,9).
\end{align*}
\end{thm}

 Guy \cite{G94} realized that $p_5+p_5+p_5$ is universal over $\Z$, and
 Sun \cite{S15} proved that the sums
\begin{align*}&p_5+p_5+2p_5,\ p_5+p_5+4p_5,\ p_5+2p_5+2p_5,
\\&p_5+2p_5+4p_5,\ p_5+p_5+5p_5,\ p_5+3p_5+6p_5
\end{align*}
are universal over $\Z$. So the converse of Theorem \ref{Th1.1}
reduces to the following conjecture of Sun.

\begin{conj}\label{Conj1.1} {\rm (Sun \cite[Remark 1.2]{S15})} The sum $p_5+bp_5+cp_5$
is universal over $\Z$ if the ordered pair $(b,c)$ is among
\begin{align*}&(1,3),\ (1,6),\ (1,8),\ (1,9),\ (1,10),\  (2,3),
\\&(2,6),\ (2,8),\ (3,3),\ (3,4),\ (3,7),\ (3,8),\ (3,9).
\end{align*}
\end{conj}

Our following result confirms this conjecture for six ordered pairs $(b,c)$ for the first time.

\begin{thm}\label{Th1.2} For
$$(b,c)=(1,3),\,(2,3),\,(2,6),\,(3,3),\,(3,4),\,(3,9),$$
the sum $p_5+bp_5+cp_5$ is universal over $\Z$.
\end{thm}

\begin{xrem} This result appeared in the initial preprint version of this paper posted to {\tt arXiv} in 2009.
\end{xrem}

Sun \cite{S15} investigated those universal sums $ap_i+bp_j+cp_k$ over
$\N$. By Sun \cite[Conjectures 1.10 and 1.13]{S15}, $p_3+p_5+p_{11}$ and
$3p_3+p_5+p_7$ should be universal over $\N$. Though we cannot prove
this, we are able to show the following result.

\begin{thm}\label{Th1.3} The sums $p_3+p_5+p_{11}$ and $3p_3+p_5+p_7$
are universal over $\Z$.
\end{thm}

Theorems \ref{Th1.2} and \ref{Th1.3} will be shown in Sections 2 and 3 respectively. Our proofs are essentially elementary and hence suitable
for general readers.

\section{Proof of Theorem \ref{Th1.2}}\setcounter{equation}{0}

\begin{lem}\label{Lem2.1} {\rm (Sun \cite[Lemma 3.2]{S15})} Let $w=x^2+3y^2\eq 4\ (\mo\ 8)$ with $x,y\in\Z$.
Then there are odd integers $u$ and $v$ such that $w=u^2+3v^2$.
\end{lem}

\begin{lem}\label{Lem2.2} Let $w=x^2+3y^2$ with $x,y$ odd and $3\nmid x$.
Then there are integers $u$ and $v$ relatively prime to $6$ such
that $w=u^2+3v^2$.
\end{lem}
\begin{proof} It suffices to consider the case $3\mid y$. Without loss of
generality, we may assume that $x\not\eq y\ (\mo\ 4)$ (otherwise we
may use $-y$ instead of $y$). Thus $(x-y)/2$ and
$(x+3y)/2=(x-y)/2+2y$ are odd. Observe that
\begin{equation}\label{2.1}x^2+3y^2=\l(\f{x+3y}2\r)^2+3\l(\f{x-y}2\r)^2.\end{equation}
As $3\nmid x$ and $3\mid y$, neither $(x-y)/2$ nor $(x+3y)/2$ is
divisible by 3. Therefore $u=(x+3y)/2$ and $v=(x-y)/2$ are
relatively prime to 6. This concludes the proof.
\end{proof}

\begin{lem}\label{Lem2.3} {\rm (Jacobi's identity)} We have
 \begin{equation}\label{2.2}3(x^2+y^2+z^2)=(x+y+z)^2+2\(\f{x+y-2z}{2}\)^2+6\(\f{x-y}{2}\)^2.\end{equation}
 \end{lem}

We need to introduce some more notation.
For $a,b,c\in\Z^+$, we set
$$E(ax^2+by^2+cz^2)=\{n\in\N:\ n\not=ax^2+by^2+cz^2\ \t{for any}\ x,y,z\in\Z\}.$$

\medskip
\noindent{\it Proof of Theorem \ref{Th1.2}}. Let
$b,c\in\Z^+$. For $n\in\N$ we have
\begin{align*} &n=p_5(x)+bp_5(y)+cp_5(z)=\f{3x^2-x}2+b\f{3y^2-y}2+c\f{3z^2-z}2
\\\iff&24n+b+c+1=(6x-1)^2+b(6y-1)^2+c(6z-1)^2.
\end{align*}
If $w\in\Z$ is relatively prime to 6, then $w$ or $-w$ is congruent to $-1$ modulo 6.
Thus, $p_5+bp_5+cp_5$ is universal over $\Z$ if and only if for any $n\in\N$ the equation
$24n+b+c+1=x^2+by^2+cz^2$ has integral solutions with $x,y,z$ relatively prime to $6$.

 Below we fix a nonnegative integer $n$.

(i) By Dickson \cite[Theorem III]{D27},
\begin{equation}\label{2.3}E(x^2+y^2+3z^2)=\{9^k(9l+6):\ k,l\in\N\}.\end{equation}
So $24n+5=u^2+v^2+3w^2$ for some $u,v,w\in\Z$.
As $3w^2\not\eq 5\ (\mo\ 4)$, $u$ or $v$ is odd.
Without loss of generality we assume that $2\nmid u$.
Since $v^2+3w^2\eq 5-u^2\eq 4\ (\mo\ 8)$, by Lemma \ref{Lem2.1}
we can rewrite $v^2+3w^2$ as $s^2+3t^2$ with $s,t$ odd.
Now we have $24n+5=u^2+s^2+3t^2$ with $u,s,t$ odd. By $u^2+s^2\eq5\eq2\ (\mo\ 3)$,
both $u$ and $s$ are relatively prime to 3. Applying Lemma \ref{Lem2.2} we can express
$s^2+3t^2$ as $y^2+3z^2$ with $y,z$ relatively prime to 6.
Thus $24n+5=u^2+y^2+3z^2$ with $u,y,z$ relatively prime to 6.
This proves the universality of $p_5+p_5+3p_5$ over $\Z$.

(ii) By Dickson \cite[Theorem X]{D27},
\begin{equation}\label{2.4}E(x^2+2y^2+3z^2)=\{4^k(16l+10):\ k,l\in\N\}.\end{equation}
So $24n+6=2u^2+v^2+3w^2$ for some $u,v,w\in\Z$. Clearly $v$ and $w$ have the same parity.
Thus $4\mid v^2+3w^2$ and hence $2u^2\eq 6\ (\mo\ 4)$. So $u$ is odd and
$v^2+3w^2\eq 6-2u^2\eq 4\ (\mo\ 8)$.
By Lemma \ref{Lem2.1} we can rewrite $v^2+3w^2$ as $s^2+3t^2$ with $s,t$ odd.
Now we have $24n+6=2u^2+s^2+3t^2$ with $u,s,t$ odd. Note that $s^2+2u^2>0$ and
$s^2+2u^2\eq0\ (\mo\ 3)$. By \cite[p.\,173]{JP} or \cite[Lemma 2.1]{S15}, we can rewrite $s^2+2u^2$ as
$x^2+2y^2$ with $x$ and $y$ relatively prime to 3. As $x^2+2y^2=s^2+2u^2\eq3\ (\mo\ 8)$,
both $x$ and $y$ are odd.
By Lemma \ref{Lem2.2}, $x^2+3t^2=r^2+3z^2$ for some integers $r,z\in\Z$ relatively prime to 6.
Thus $24n+6=r^2+2y^2+3z^2$ with $r,y,z$ relatively prime to 6. It follows that $p_5+2p_5+3p_5$ is universal over $\Z$.

(iii) By Dickson \cite[Theorem IV]{D27},
\begin{equation}\label{2.5}E(x^2+3y^2+3z^2)=\{9^k(3l+2):\ k,l\in\N\}.\end{equation}
So $24n+7=u^2+3v^2+3w^2$ for some $u,v,w\in\Z$. Since $u^2\not\eq 7\ (\mo\ 4)$,
without loss of generality we assume that $2\nmid w$. As $u^2+3v^2\eq 7-3w^2\eq 4\ (\mo\ 8)$,
by Lemma \ref{Lem2.1} there are odd integers $s$ and $t$ such that $u^2+3v^2=s^2+3t^2$.
Thus $24n+7=s^2+3t^2+3w^2$ with $s,t,w$ odd. Clearly, $s$ is relatively prime to 6.
By Lemma \ref{Lem2.2}, $s^2+3t^2=x_0^2+3y^2$ for some integers $x_0$ and $y$ relatively prime to 6,
and $x_0^2+3w^2=x^2+3z^2$ for some integers $x$ and $z$ relatively prime to 6. Therefore
$24n+7=x^2+3y^2+3z^2$ with $x,y,z$ relatively prime to 6. This proves the universality of
$p_5+3p_5+3p_5$ over $\Z$.

(iv) By \cite[Theorem 1.7(iii)]{S15}, $24n+8=u^2+v^2+3w^2$ for some $u,v,w\in\Z$ with $2\nmid w$.
Clearly $u\not\eq v\ (\mo\ 2)$. Without loss of generality, we assume that $u=2r$ with $r\in\Z$.
Since $(2r)^2+v^2\eq 8\eq2\ (\mo\ 3)$, both $r$ and $v$ are relatively prime to 3.
As $v$ and $w$ are odd, $v^2+3w^2\eq 4\ (\mo\ 8)$ and hence $r$ is odd. By Lemma \ref{Lem2.2}, we can rewrite
$v^2+3w^2$ as $x^2+3y^2$ with $x$ and $y$ relatively prime to 6. Note that
$24n+8=4r^2+v^2+3w^2=x^2+3y^2+4r^2$ with $x,y,r$ relatively prime to 6.
It follows that $p_5+3p_5+4p_5$ is universal over $\Z$.

(v) By (2.3), $24n+13=u^2+v^2+3w^2$ for some $u,v,w\in\Z$. Since $3w^2\not\eq 13\eq1\ (\mo\ 4)$,
without loss of generality we may assume that $u$ is odd. As $v^2+3w^2\eq 13-u^2\eq 4\ (\mo\ 8)$,
by Lemma \ref{Lem2.1} we can rewrite $v^2+3w^2$ as $s^2+3t^2$ with $s$ and $t$ odd. Thus
$24n+13=u^2+s^2+3t^2$ with $u,s,t$ odd. Since $u^2+s^2\eq 13\eq1\ (\mo\ 3)$,
without loss of generality we may assume that $3\nmid u$ and $s=3r$ with $r\in\Z$.
By Lemma \ref{Lem2.2}, $u^2+3t^2=x^2+3y_0^2$ for some integers $x$ and $y_0$ relatively prime to 6,
also $y_0^2+3r^2=y^2+3z^2$ for some integers $y$ and $z$ relatively prime to 6. Thus
$24n+13=x^2+3y_0^2+9r^2=x^2+3y^2+9z^2$ with $x,y,z$ relatively prime to 6.
This proves the universality of $p_5+3p_5+9p_5$ over $\Z$.

(vi) By the Gauss-Legendre theorem (cf. \cite[pp.\,17-23]{N96}),
$8n+3=x^2+y^2+z^2$ for some odd integers $x,y,z$. Without loss of
generality we may assume that $x\not\eq y\ (\mo\ 4)$. By Jacobi's
identity (\ref{2.2}), we have $3(8n+3)=u^2+2v^2+6w^2$, where $u=x+y+z$,
$v=(x+y)/2-z$ and $w=(x-y)/2$ are odd integers. As $u^2+2v^2$ is a
positive integer divisible by 3, by \cite[p.\,173]{JP} or \cite[Lemma 2.1]{S15} we can write
$u^2+2v^2=a^2+2b^2$ with $a$ and $b$ relatively prime to 3. Since
$a^2+2b^2=u^2+2v^2\eq 3\ (\mo\ 8)$, both $a$ and $b$ are odd. By
Lemma \ref{Lem2.2}, $b^2+3w^2=c^2+3d^2$ for some integers $c$ and $d$
relatively prime to 6. Thus $24n+9=a^2+2b^2+6w^2=a^2+2c^2+6d^2$ with
$a,c,d$ relatively prime to 6. It follows that $p_5+2p_5+6p_5$ is
universal over $\Z$.

In view of the above, we have completed the proof of Theorem \ref{Th1.2}.
\qed

\section{Proof of Theorem \ref{Th1.3}}\setcounter{equation}{0}

\medskip
\noindent{\it Proof of Theorem \ref{Th1.3}}. (i) Let $n\in\N$. By part (v) in the proof of Theorem \ref{Th1.2},
there are integers $u,v,w\in\Z$ relatively prime to 6 such that
$$72n+61=24(3n+2)+13=9u^2+3v^2+w^2.$$
Clearly $w^2\eq 61-3v^2\eq 7^2\ (\mo\ 9)$ and hence $w\eq\pm 7\ (\mo\ 9)$. So there are $x,y,z\in\Z$ such that
$$72n+61=9(2x+1)^2+3(6y-1)^2+(18z-7)^2$$
and hence $n=p_3(x)+p_5(y)+p_{11}(z)$. (Note that
$p_{11}(x)=9(x^2-x)/2+x=(9x^2-7x)/2$.)

(ii) Let $n\in\N$. It is easy to see that
\begin{align*} &n=3p_3(x)+p_5(y)+p_7(z)
\\\iff&120n+77=5(3(2x+1))^2+5(6y-1)^2+3(10z-3)^2.
\end{align*}

Suppose $120n+77=5x^2+5y^2+3z^2$ for some $x,y,z\in\Z$ with $z$ odd. Then
$x^2+y^2\eq 77-3z^2\eq 2\ (\mo\ 4)$ and hence $x$ and $y$ are odd.
Note that $3z^2\eq 77\eq 12\ (\mo\ 5)$ and hence $z\eq\pm3\ (\mo\ 10)$.
As $5x^2+5y^2\eq 77\eq 5\ (\mo\ 3)$, exactly one of $x$ and $y$ is divisible by 3.
Thus there are $u,v,w\in\Z$ such that
$$120n+77=5(3(2u+1))^2+5(6v-1)^2+3(10w-3)^2.$$

 By the above, to prove the universality of $3p_3+p_5+p_7$ over $\Z$,
 we only need to show that $120n+77=5x^2+5y^2+3z^2$ for some $x,y,z\in\Z$ with $z$ odd.

By (\ref{2.3}), there are $u,v,w\in\Z$ such that $120n+77=u^2+v^2+3w^2$. As $3w^2\not\eq 77\eq1\ (\mo\ 4)$,
$u$ or $v$ is odd, say, $2\nmid u$. As $v^2+3w^2\eq 77-u^2\eq 4\ (\mo\ 8)$,
by Lemma \ref{Lem2.1} we may assume that $v$ and $w$ are odd without loss of generality.

We claim that $120n+77=a^2+b^2+3c^2$ for some odd integers $a,b,c$ with $c\eq\pm 2\ (\mo\ 5)$.
This holds if $w\eq\pm 2\ (\mo\ 5)$. Suppose that $w\not\eq\pm2\ (\mo\ 5)$.
If $w\eq\pm1\ (\mo\ 5)$, then $u^2+v^2\eq 77-3w^2\eq-1\ (\mo\ 5)$ and hence
$u$ or $v$ is divisible by 5. If $w\eq0\ (\mo\ 5)$, then $u^2+v^2\eq 77\eq2\ (\mo\ 5)$
and hence $u^2\eq v^2\eq1\ (\mo\ 5)$.
Without loss of generality, we assume that one of $v$ and $w$ is divisible by 5
and the other one is congruent to $1$ or $-1$ modulo 5, we may also suppose that $v\not\eq w\ (\mo\ 4)$
(otherwise we may use $-w$ instead of $w$). By the identity (\ref{2.1}),
$$v^2+3w^2=\l(\f{v+3w}2\r)^2+3\l(\f{v-w}2\r)^2.$$
Note that both $(v-w)/2$ and $(v+3w)/2=(v-w)/2+2w$ are odd. Also, $(v-w)/2$ is congruent to $2$ or $-2$ modulo 5.
This confirms the claim.

By the above, there are odd integers $a,b,c\in\Z$ with $c\eq\pm 2\ (\mo\ 5)$ such that
$120n+77=a^2+b^2+3c^2$. Since $3c^2\eq 77\ (\mo\ 5)$, we have $5\mid a^2+b^2$ and hence
$a^2\eq(2b)^2\ (\mo\ 5)$. Without loss of generality we assume that $a\eq 2b\ (\mo\ 5)$.
Then $x=(2a+b)/5$ and $y=(a-2b)/5$ are odd integers, and
$$a^2+b^2=(2x+y)^2+(x-2y)^2=5(x^2+y^2).$$
Now we have $120n+77=5(x^2+y)^2+3c^2$ with $x,y,c$ odd.

This concludes our proof of Theorem 1.4. \qed

\subsection*{Acknowledgements}
The authors would like to thank Dr. Hao Pan for his helpful comments.
The second author was supported by the National Natural Science Foundation of China (grant 11571162).


\medskip

\begin{thebibliography}{99}


\bibitem{CH} W. K. Chan and A. Haensch, {\it Alomost universal ternary sums of squares and triangular numbers}, in: Quadratic and Higher Degree Forms,  Dev. Math., Vol. 31, Springer, New York, 2013, pp. 51--62.

\bibitem{CO} W. K. Chan and B.-K. Oh, {\it Almost universal ternary sums of triangular numbers}, Proc. Amer. Math. Soc. {\bf 137} (2009), 3553--3562.

\bibitem{D27} L. E. Dickson, {\it Integers represented by positive ternary quadratic forms}, Bull. Amer. Math. Soc. {\bf 33} (1927), 63--77.

\bibitem{D99} L. E. Dickson, History of the Theory of Numbers, {\rm Vol. II}, AMS Chelsea Publ., 1999.

\bibitem{GPS} S. Guo, H. Pan and Z.-W. Sun, {\it Mixed sums of squares and triangular numbers (II)}, Integers {\bf  7} (2007), \#A56, 5pp (electronic).

\bibitem{G94} R. K. Guy, {\it Every number is expressible as the sum of how many polygonal numbers?} Amer. Math. Monthly {\bf 101} (1994), 169--172.

\bibitem{JP} B. W. Jones and G. Pall, {\it Regular and semi-regular positive ternary quadratic forms}, Acta Math. {\bf 70} (1939), 165--191.

\bibitem{KS} B. Kane and Z.-W. Sun, {\it On almost universal mixed sums of squares and triangular numbers}, Trans. Amer. Math. Soc. {\bf 362} (2010), 6425--6455.

\bibitem{OS} B.-K. Oh and Z.-W. Sun, {\it Mixed sums of squares and triangular numbers (III)}, J. Number Theory {\bf 129} (2009), 964-969.

\bibitem{MW} C. J. Moreno and S. S. Wagstaff, Sums of Squares of Integers, Chapman \& Hall/CRC, Boca Raton, FL, 2006.

\bibitem{N96} M. B. Nathanson, Additive Number Theory: The
Classical Bases, Grad. Texts in Math., vol. 164, Springer,
New York, 1996.

\bibitem{S07} Z.-W. Sun, {\it Mixed sums of
squares and triangular numbers}, Acta Arith. {\bf 127} (2007), 103--113.

\bibitem{S15} Z.-W. Sun, {\it On universal sums of polygonal numbers},
Sci. China Math. {\bf 58} (2015), 1367--1396.

\end{thebibliography}
\end{document}